\documentclass[letterpaper]{article}
\usepackage{amsmath,amsfonts,amsthm,xypic,hyperref}

{}
\newtheorem{theorem}{Theorem}[section]
\newtheorem{corollary}[theorem]{Corollary}

\newtheorem{prop}[]{Proposition}

\newcommand{\nc}{\newcommand}
\nc{\FH}{\mathcal H} \nc{\CC}{\mathbb C} \nc{\JJ}{\mathcal J}
\nc{\KK}{\mathbb K} \nc{\RR}{\mathbb R} \nc{\LL}{\mathcal L}
\nc{\Ll}{\ell} \nc{\NN}{\mathbb N} \nc{\ZZ}{\mathbb Z} \nc
{\HH}{\mathbb H} \nc {\OO}{\mathcal{O}} \nc{\lra}{\longrightarrow}
\nc{\bdot}{\bullet} \nc{\w}{\omega} \nc{\dd}{\mathcal{D}}
\newcommand{\vol}{\mathrm{vol}}

\newcommand{\End}{\mathrm{End}}

\newcommand{\delbar}{\overline{\partial}}
\newcommand{\so}{\mathfrak{so}}
\newcommand{\IP}[1]{\langle #1\rangle}
\nc{\ba}{\overline} \nc{\del}{\partial}
\nc{\de}{\delta}\nc{\debar}{\overline{\delta}}
\nc{\Ja}{e^{\tfrac{\pi}{2}\JJ_1}}
\nc{\Jb}{e^{\tfrac{\pi}{2}\JJ_2}}

\author{Marco Gualtieri\\ \textsl{mgualtie@fields.utoronto.ca}}
\title{\bf Generalized geometry and the Hodge decomposition}
\begin{document}
\date{}
\maketitle
\begin{abstract}
In this lecture, delivered at the string theory and geometry
workshop in Oberwolfach, we review some of the concepts of
generalized geometry, as introduced by Hitchin and developed in
the speaker's thesis. We also prove a Hodge decomposition for the
twisted cohomology of a compact generalized K\"ahler manifold, as
well as a generalization of the $dd^c$-lemma of K\"ahler geometry.
\end{abstract}
\section{Geometry of $T\oplus T^*$}
The sum $T\oplus T^*$ of the tangent and cotangent bundles of an
$n$-dimensional manifold has a natural $O(n,n)$ structure given by
the inner product
\[
\IP{X+\xi,Y+\eta}=\tfrac{1}{2}(\xi(Y)+\eta(X)),
\]
and we may describe the Lie algebra in the usual way:
\[
\mathfrak{so}(n,n)=\wedge^2 T\oplus \End(T)\oplus \wedge^2 T^*.
\]
Hence we may view 2-forms $B$ and bivectors $\beta$ as
infinitesimal symmetries of $T\oplus T^*$.  We may also form the
Clifford algebra $CL(T\oplus T^*)$, which has a spin
representation on the Clifford module $\wedge^\bullet T^*$ as
described by Cartan:
\[
(X+\xi)\cdot\rho = i_X\rho + \xi\wedge \rho,
\]
for $X+\xi\in T\oplus T^*$ and $\rho\in\wedge^\bullet T^*$.  This
means that we may view differential forms as
spinors\footnote{Actually, the bundle of spinors differs from
$\wedge^\bullet T^*$ by tensoring with the line bundle $\det
T^{1/2}$; we assume a trivialization has been chosen -- this is
related to the physicists' \emph{dilaton} field.} for $T\oplus
T^*$. From the general theory of spinors, this implies that there
is a $Spin_{o}(n,n)$-invariant bilinear form
\[
\IP{,}\colon\wedge^\bullet T^* \times \wedge^\bullet T^* \lra \det
T^*,
\]
given by $\IP{\alpha,\beta} = [\alpha\wedge\sigma(\beta)]_n$,
where $\sigma$ is the anti-automorphism which reverses the order
of wedge product.

Another structure emerging from the interpretation of forms as
spinors is the Courant bracket $[,]_H$ on sections of $T\oplus
T^*$, obtained as the derived bracket (see~\cite{KS}) of the
natural differential operator $d+H\wedge\cdot$ acting on
differential forms, where $d$ is the exterior derivative and $H\in
\Omega^3_{cl}(M)$.  When $H=0$, we have the following:
\begin{prop}The group of orthogonal automorphisms of the Courant bracket for
$H=0$ is a semidirect product of $\mathrm{Diff}(M)$ and
$\Omega^2_{cl}(M)$, where $B\in\Omega^2_{cl}(M)$ acts as the shear
$\exp(B)$ on $T\oplus T^*$.
\end{prop}

In this way we see that the natural appearance of
$H\in\Omega^3_{cl}(M)$ and the action of $B\in\Omega^2_{cl}(M)$
coincide precisely with the physicists' description of the
Neveu--Schwarz 3-form flux and the action of the B-field,
respectively.

\section{Generalized complex geometry}

A generalized complex structure is an integrable reduction of the
structure group of $T\oplus T^*$ from $O(2n,2n)$ to $U(n,n)$ (only
possible when $\dim_\RR M = 2n$).  This is equivalent to the
choice of an orthogonal complex structure
\[
\JJ\in O(T\oplus T^*),\ \JJ^2=-1.
\]
The integrability condition is that the $+i$-eigenbundle of $\JJ$,
\[
E<(T\oplus T^*)\otimes \CC,
\]
must be closed under the Courant bracket.  If $H$ is nonzero we
call this a twisted generalized complex structure.  The Courant
bracket is a Lie bracket when restricted to $E$ and therefore we
may form the associated differential graded algebra:
\[
\mathcal{E}=(\wedge^\bullet E^*, d_E).
\]
Because $E^*$ is identified with $\overline{E}$ by the metric on
$T\oplus T^*$, we see that it also acquires a Lie bracket. By a
general result of Lu, Weinstein, and Xu \cite{Lu}, this Lie
bracket makes $\mathcal{E}$ into a differential Gerstenhaber
algebra.
\begin{theorem}
The differential Gerstenhaber algebra $\mathcal{E}$ is elliptic
and it gives rise to a Kuranishi deformation theory for any
generalized complex structure. The tangent space to the
deformation space, in the unobstructed case, is
$H^2(\mathcal{E})$, and obstructions lie in $H^3(\mathcal{E})$.
\end{theorem}

For example, let $J\in\End(T)$ be a usual complex structure, and
form the generalized complex structure
\[
\JJ=\begin{pmatrix}-J&0\\0&J^*\end{pmatrix}.
\]
Then, $E=T_{0,1}\oplus T^*_{1,0}$, so that $\mathcal{E}$ is simply
the Dolbeault complex of the holomorphic multivectors.
Consequently,
\[
H^2(\mathcal{E})=H^0(\wedge^2 T)\oplus H^1(T)\oplus
H^2(\mathcal{O}).
\]
For a complex surface, a holomorphic bivector $\beta$ always
integrates to an actual deformation, and so for $\CC P^2$, for
example, we obtain a new generalized complex structure which is
complex along an anticanonical divisor (the vanishing locus of
$\beta$) and the B-field transform of a symplectic structure away
from the cubic.

This provides an alternative interpretation of the extended
deformation parameter $\beta$, which is normally viewed as a
noncommutative deformation of the algebra defining $\CC P^2$. Note
that the usual translation parameter along the commutative
elliptic curve can be obtained by differentiating $\beta$ along
its vanishing set.

The previous example indicates that the algebraic type of a
generalized complex structure may jump along loci in the manifold.
Indeed, a generalized complex structure on a $2n$-manifold may
have types ${0,\cdots, n}$, with $0$ denoting the (generic)
symplectic type and $n$ denoting the complex type.  Type may jump
up, but only by an even number.

\begin{theorem}[Generalized Darboux theorem~\cite{Gu}]
Away from type jumping loci, a generalized complex manifold of
type $k$ is locally isomorphic, via a diffeomorphism and a
$B$-field transform, to $\CC^k\times \RR^{2n-2k}_{\omega_0}$,
where $\omega_0$ is the usual Darboux symplectic form.
\end{theorem}

Generalized complex manifolds also have natural sub-objects,
called generalized complex submanifolds~\cite{Gu}.  These
sub-objects correspond exactly with the physicists' notion of
topological D-branes; in particular, one recovers, in the
symplectic case, the co-isotropic A-branes of Kapustin and
Orlov~\cite{Ka}. There is also a natural notion of generalized
holomorphic bundle supported on a generalized complex submanifold,
a concept which seems to correspond to D-branes of higher rank.
One can even see how such a brane could deform into several branes
of lower rank.

\section{Generalized Riemannian geometry and the Born-Infeld metric}

A generalized Riemannian metric is a reduction of the structure
group of $T\oplus T^*$ from $O(n,n)$ to $O(n)\times O(n)$.  This
is equivalent to specifying a maximal positive-definite subbundle,
$C_+<T\oplus T^*$, which can be described as the graph of $b+g$,
where $g$ is a usual Riemannian metric and $b$ is a 2-form, each
viewed as defining a map $T\lra T^*$ via interior product.  The
graph of $b-g$ is denoted by $C_-$, the orthogonal complement to
$C_+$.  These data determine a positive-definite metric on
$T\oplus T^*$ by simply taking $\IP{,}|_{C_+} - \IP{,}|_{C_-}$.
This metric, evaluated on $A,B\in T\oplus T^*$, can be written as
$\IP{GA,B}$, where $G$ is the obvious involution, expressible in
terms of the data as follows:
\[
G = \left(
\begin{array}{cc}
  1 &  \\
  b & 1 \\
\end{array}
\right) \left(
\begin{array}{cc}
   & g^{-1}  \\
  g &  \\
\end{array}\right)\left(
\begin{array}{cc}
  1 &  \\
-b & 1 \\
\end{array}%
\right) = \left(%
\begin{array}{cc}
  -g^{-1}b & g^{-1}\\
  g-bg^{-1}b & bg^{-1} \\
\end{array}
\right).
\]
The restriction of this metric to the subbundle $T$ is the
Riemannian metric $g-bg^{-1}b$.  Note that the volume form induced
by this last metric is
\[
\vol_G=\det(g-bg^{-1}b)^{1/2} = (\det(g)
\det(1-g^{-1}bg^{-1}b))^{1/2} = \frac{\det(g+b)}{{\det g}^{1/2}}.
\]

Let $\ast=a_1\cdots a_n$ be the product in $CL(C_+)<CL(T\oplus
T^*)$ of an oriented orthonormal basis for $C_+$.  This volume
element acts on the differential forms via the spin
representation, and is related to the Hodge star operator $\star$
in the following way: if $b=0$ then
\[
\star\rho = \sigma(\sigma(\ast)\cdot\rho).
\]
Since $\ast^2 = (-1)^{n(n-1)/2}$ and $\IP{\alpha,\beta} =
(-1)^{n(n-1)/2}\IP{\beta,\alpha}$, we see that the volume form,
\[
\IP{\alpha,\sigma(\ast)\beta},
\]
is symmetric in $\alpha,\beta$.  For $b=0$, it is nothing but the
Hodge volume:
\[
\IP{\alpha,\sigma(\ast)\beta}=\alpha\wedge\star\beta =
g(\alpha,\beta)\vol_g,
\]
where $g(\alpha,\beta)$ is the positive-definite Hodge metric on
differential forms.  In the general case, we obtain a different
symmetric positive-definite volume form,
\[
\IP{\alpha,\sigma(\ast)\beta} = G(\alpha,\beta)\IP{1,\sigma(\ast)
1}=G(\alpha,\beta)\frac{\det(g+b)}{{\det g}^{1/2}} =
G(\alpha,\beta)\vol_G,
\]
where $G(\alpha,\beta)$ is a positive-definite metric on forms
satisfying $(1,1)=1$.  We call this expression the Born-Infeld
volume, to coincide with physics terminology.  Therefore, for any
generalized Riemannian structure, we may define the following
positive-definite Hermitian inner product on differential forms:
\[
h(\alpha,\beta) = \int_M\IP{\alpha,\sigma(\ast)\bar\beta},
\]
which we call the Born-Infeld inner product.  It is a direct
generalization of the Hodge inner product of Riemannian geometry.

To develop Hodge theory for generalized Riemannian manifolds, we
calculate the adjoint of the twisted exterior derivative $d_H$.
Note first that the exterior derivative is such that
\[
\IP{d_H\alpha,\beta}- (-1)^{\dim_\RR M}\IP{\alpha,d_H\beta}\ \
\text{is exact,}
\]
so that for a compact even-dimensional manifold,
\[
\int_M\IP{d_H\alpha,\beta} = \int_M\IP{\alpha,d_H\beta}.
\]

With this in mind, we may determine the Born-Infeld adjoint of
$d_H$:
\begin{align*}
h(d_H\alpha,\beta)&=\int_M\IP{d_H\alpha, \sigma(\ast)\bar\beta}\\
&=\int_M\IP{\alpha, d_H\sigma(\ast)\bar\beta}\\
&=\int_M\IP{\alpha, \sigma(\ast)\ast d_H\sigma(\ast)\bar\beta}\\
&=h(\alpha,\ast d_H\sigma(\ast)\beta),
\end{align*}
proving that, for an even-dimensional manifold,
\[
d_H^\ast = \ast\cdot d_H\cdot\sigma(\ast) = \ast\cdot d_H\cdot
\ast^{-1}.
\]

As in the Riemannian case, the operator $\mathcal{D}_+= d_H+d_H^*$
is elliptic (as an operator
$\wedge^{ev/od}T^*\rightarrow\wedge^{od/ev}T^*$) and so,
therefore, is the Laplacian
\[
\Delta_{d_H}=\mathcal{D}_+^2 = d_Hd^*_H+d^*_Hd_H.
\]
Proceeding in the usual way, we may conclude that on a compact
generalized Riemannian manifold, every $H$-twisted de Rham
cohomology class has a unique $\Delta_{d_H}$-harmonic
representative.  There is a gauge freedom here, in the sense that,
given any 2-form $b'$, the automorphism $e^{b'}$ takes harmonic
representatives for $(g,b,H)$ to those for $(g,b+b', H-db')$.

\section{Generalized K\"ahler geometry and the Hodge decomposition}

A generalized K\"ahler structure is a further integrable reduction
of the structure group of $T\oplus T^*$ to $U(n)\times U(n)$. As
defined in~\cite{Gu}, it consists of two commuting generalized
complex structures $(\JJ_1,\JJ_2)$ such that the involution
$-\JJ_1\JJ_2 = G$ determines a generalized Riemannian metric on
$T\oplus T^*$. The standard example of such a pair is obtained
from a usual K\"ahler structure $(g,J,\omega)$, i.e. a Riemannian
metric $g$, a complex structure $J$, and a symplectic structure
$\omega$, such that the following diagram commutes:
\begin{equation*}
\xymatrix{T\ar[rr]^{g}&&T^*\\
&T\ar[ul]^{J} \ar[ur]_{\omega}&}
\end{equation*}
By taking
\[
\JJ_1=\left(%
\begin{array}{cc}
  -J &  \\
   & J^* \\
\end{array}
\right),\ \ \ \ \JJ_2 = \left(%
\begin{array}{cc}
    & \omega^{-1} \\
  -\omega &  \\
\end{array}%
\right),
\]
we see that these generalized complex structures commute and
\[
-\JJ_1\JJ_2 = \left(%
\begin{array}{cc}
    & g^{-1} \\
  g &  \\
\end{array}%
\right),
\]
defining a generalized Riemannian metric with $b=0$.

In the preceding example, the types of the generalized complex
structures $(\JJ_1,\JJ_2)$ are $(n,0)$, since one is complex and
the other is symplectic.  In general, though type jumping may
occur, we have the following constraints on the pair of types:
\begin{align*}
&\text{type}(\JJ_1)+\text{type}(\JJ_2) \equiv n \pmod 2,\
\text{and}\\
&\text{type}(\JJ_1)+\text{type}(\JJ_2) \leq n. \end{align*}

In~\cite{Gu}, it is proven that generalized K\"ahler geometry is
equivalent to a bi-Hermitian geometry with torsion, first
described by Gates, Hull, and Ro\v{c}ek~\cite{Ga} in their study
of $N=(2,2)$ supersymmetric sigma models.  This equivalence
indicates how it is possible to `topologically twist' these models
in general.  In what follows, we are interested in what
implications generalized K\"ahler geometry has for differential
forms, and in particular whether there is a generalization of the
Hodge decomposition.

First observe that both $\JJ_1,\JJ_2$ are in $\so(T\oplus T^*)$,
and via the spin representation they act on differential forms.
$\JJ_1$ induces a decomposition of forms into its eigenspaces:
\[
\wedge^\bullet T^*\otimes\CC = U_{-n}\oplus\cdots\oplus
U_0\oplus\cdots\oplus U_n,
\]
where $U_k$ is the $ik$-eigenspace of $\JJ_1$.  Furthermore, the
exterior derivative $d_H$, acting on sections of $U_k$, decomposes
as the sum of the two projections $\delbar_1,\del_1$, to
$U_{k+1},U_{k-1}$, respectively.  That is,
\begin{equation*}
\xymatrix{C^\infty(U_k)\ar@<0.5ex>[r]^{\delbar_1}&C^\infty(U_{k+1})\ar@<0.5ex>[l]^{\del_1}}.
\end{equation*}
The commuting endomorphism $\JJ_2$ engenders a further
decomposition of the $U_k$:
\[
U_k = U_{k,|k|-n}\oplus U_{k,|k|-n+2}\oplus \cdots \oplus
U_{k,n-|k|},
\]
where $U_{p,q}$ is the intersection of the $ip$-eigenspace of
$\JJ_1$ and the $iq$-eigenspace of $\JJ_2$. In this way we obtain
a $(p,q)$ decomposition of the differential forms into the
following diamond:
\[
\begin{array}{ccccccc}
  &  &  &  U_{0,n}&  &  &  \\
   &  & \cdots &  & \cdots &  &  \\
   &U_{-n+1,1}&  &  &  &U_{n-1,1}  &  \\
  U_{-n,0} &  &  & \cdots &  &  & U_{n,0} \\
   & U_{-n+1,-1} &  &  &  &U_{n-1,-1}  &  \\
   &  & \cdots & & \cdots &  &  \\
   &  &  &U_{0,-n}  &  &  &  \\
\end{array}
\]
This decomposition is orthogonal with respect to the Born-Infeld
metric, and gives rise to the following decomposition of the
exterior derivative:
\[
d_H = \de_++\de_- + \debar_+ + \debar_-,
\]
where the differential operators act as follows:
\[
\xymatrix{
U_{p-1,q+1} & &U_{p+1,q+1} \\
 &U_{p,q}\ar@{.>}[r]^{\delbar_1}\ar@{.>}[l]^{\del_1}\ar@{.>}[u]^{\delbar_2}\ar@{.>}[d]^{\del_2}\ar[lu]^{\de_-}\ar[ru]^{\debar_+}\ar[ld]^{\de_+}\ar[rd]^{\debar_-} & \\
 U_{p-1,q-1}& & U_{p+1,q-1}
}
\]
and where we have, for definiteness, $\delbar_1= \debar_+ +
\debar_-$ and $\delbar_2 = \debar_+ + \de_-$.

The following proposition gives the crucial relationship between
these operators, and is a generalization of the usual K\"ahler
identities:
\begin{prop}[generalized K\"ahler identities~\cite{Gu2}]
For a generalized K\"ahler structure, we have the identities
\[
\debar_+^* = -\de_+\ \ \text{and}\ \ \debar_-^* = \de_-.
\]
\end{prop}
\noindent These simple identities imply the equality of all
available Laplacians:
\[
\Delta_{d_H}=2\Delta_{\delbar_{1/2}}=2\Delta_{\del_{1/2}}=4\Delta_{\debar_{\pm}}=4\Delta_{\de_{\pm}},
\]
and so, finally, we obtain a $(p,q)$ decomposition for the twisted
cohomology of any compact generalized K\"ahler manifold.
\begin{theorem}[Hodge decomposition~\cite{Gu2}]
The twisted cohomology of a compact twisted $2n$-dimensional
generalized K\"ahler manifold carries a Hodge decomposition:
\[
H^\bullet_H(M,\CC) = \bigoplus_{\substack{|p+q|\leq n\\
p+q\equiv n(\text{mod}\ 2)}} \mathcal{H}^{p,q},
\]
where $\mathcal{H}^{p,q}$ are $\Delta_{d_H}$-harmonic forms in
$U_{p,q}$.
\end{theorem}
Note that in the usual K\"ahler case, this $(p,q)$ decomposition
is \emph{not} the Dolbeault decomposition: it was called the
Clifford decomposition by Michelsohn~\cite{Mi}, and there is an
orthogonal transformation, called the Hodge automorphism, taking
it to the usual Dolbeault decomposition.  A striking feature of
the Clifford decomposition is that a form of type $(p,q)$ is
closed if and only if it is co-closed and hence harmonic.

A consequence of this result is a generalization of the well-known
$\del\delbar$-lemma of K\"ahler geometry.  Any (twisted)
generalized complex structure gives rise to a real differential
operator $d^\JJ_H = [d_H,\JJ]$ on forms.  In the complex case
$d^\JJ = d^c= i(\delbar-\del)$, whereas in the symplectic case,
$d^\JJ = \delta$, the Koszul symplectic adjoint of $d$.  The
$dd^\JJ$-property, studied in more detail by Cavalcanti~\cite{Ca},
is then defined as follows.
\begin{paragraph}{Definition
($dd^\JJ$ property):} {\it A generalized complex manifold
$(M,\JJ)$ satisfies the $dd^J$ property iff the following are
equivalent:
\begin{itemize}
\item $\rho$ is $d$-closed and $d^\JJ$-exact,
\item $\rho$ is $d^\JJ$-closed and $d$-exact,
\item $\rho = dd^\JJ\tau$ for some $\tau$.\nopagebreak
\end{itemize}}
\end{paragraph}
We have given the property for $H=0$; in general simply replace
$d$ by $d_H$.  Now we can state the first corollary of the
previous theorem:
\begin{corollary}[$dd^\JJ$ lemma~\cite{Gu2}]
A compact twisted generalized K\"ahler manifold satisfies the
$d_Hd^\JJ_H$ property with respect to both $\JJ_1$ and $\JJ_2$.
\end{corollary}
In the usual K\"ahler case, this implies that both the $dd^c$ and
$d\delta$-lemmas are satisfied.  As shown by Merkulov
(see~\cite{Ca}), the $d\delta$-lemma is equivalent to the strong
Lefschetz property, which we know is satisfied by any compact
K\"ahler manifold.

A second corollary of the Hodge decomposition concerns a
generalization of the fact that the odd Betti numbers of a compact
K\"ahler manifold must be even.  Observing that
$\overline{\mathcal{H}^{p,q}} = \mathcal{H}^{-p,-q}$, we obtain a
constraint on the parity of the even or odd twisted Betti numbers
$b^{ev/od}_H = \dim H^{ev/od}_H(M)$:
\begin{corollary}
Let $M$ be a compact twisted generalized K\"ahler manifold.  If
$\dim M=4k+2$, then both $b^{od}_H$ and $b^{ev}_H$ must be even.
If $\dim_\RR M = 4k$, then the generalized K\"ahler pair may have
types of parity either $(od,od)$ or $(ev,ev)$.  In the former
case, $b^{ev}_H$ must be even, whereas in the latter case,
$b^{od}_H$ must be even.
\end{corollary}
By applying this corollary, we see at once that the 4-manifold
$\CC P^2$ does not admit a generalized K\"ahler structure with
types $(1,1)$.

\end{document}